\documentclass[10pt,twocolumn,twoside]{IEEEtran}

\usepackage{bm}   
\usepackage{amsmath}
\allowdisplaybreaks[3]
\usepackage{amssymb}
\usepackage{indentfirst}  
\usepackage{graphicx}  
\usepackage{subfigure}
\usepackage{multirow} 
\usepackage{booktabs} 
\usepackage{cite}  
\usepackage{threeparttable}   
\usepackage{algorithm}
\usepackage{algpseudocode}
\usepackage{arydshln}
\usepackage{color}
\usepackage[american]{babel}
\usepackage{microtype}
\usepackage[numbers,sort&compress]{natbib}  
\usepackage{mathtools}
\usepackage{epstopdf}

\newcommand{\tabincell}[2]{\begin{tabular}{@{}#1@{}}#2\end{tabular}}  

\newtheorem{remark}{Remark}

\newtheorem{prop}{Proposition}

\def\mV{\mathcal{V}}
\def\mE{\mathcal{E}}

\def\mN{\mathcal{N}}

\def\mbR{\mathbb{R}}
\def\mbC{\mathbb{C}}

\def\mEf{\mE_{\textup{flex}}}

\def\upj{\textup{j}}
\def\opt{\textup{opt}}
\def\sdp{\textup{sdp}}

\usepackage[font=footnotesize,justification=centering]{caption}
\begin{document}
%
\title{Convex Relaxation of AC Optimal Power Flow with Flexible Transmission Line Impedances}

\author{Yue~Song,~\IEEEmembership{Member,~IEEE,}
        David~J.~Hill,~\IEEEmembership{Life Fellow,~IEEE,}
        Tao~Liu,~\IEEEmembership{Member,~IEEE,}\\
        and~Tianlun~Chen,~\IEEEmembership{Member,~IEEE}

\vspace{-10pt}

\thanks{This work was supported by the HKU Seed Fund for Basic Research for New Staff under Project No. 202009185007.}
\thanks{Y. Song, T. Liu, and T. Chen are with the Department of Electrical and Electronic Engineering, The University of Hong Kong, Hong Kong (e-mail: yuesong@eee.hku.hk; taoliu@eee.hku.hk; tlchen@eee.hku.hk).}
\thanks{D. J. Hill is with the Department of Electrical and Electronic Engineering, The University of Hong Kong, Hong Kong, and also with the School of Electrical Engineering and Telecommunications, The University of New South Wales, Kensington, NSW 2052, Australia (e-mail: dhill@eee.hku.hk).}
}

\markboth{IEEE TRANSACTIONS ON POWER SYSTEMS} 
{Song \MakeLowercase{\textit{et al.}}: Convex Relaxation of AC OPF with Flexible Line Impedances}

\maketitle

\begin{abstract}
  Flexible transmission line impedances on one hand are a promising control resource for facilitating grid flexibility, but on the other hand add much complexity to the concerned optimization problems.
  This paper develops a convexification method for the AC optimal power flow with flexible line impedances.
  First, it is discovered that a flexible-impedance line is equivalent to a constant-impedance line linking a pair of transformers with correlated and continuously adjustable tap ratios.
  Then, with this circuit equivalent, the original optimization problem is reformulated into a semi-definite program under the existing convex relaxation framework, which improves the solution tractability and optimality in an easy-to-implement manner.
  The proposed method is verified by numerical tests on the IEEE 118-bus system.
\end{abstract}


%
\IEEEpeerreviewmaketitle


\section{Introduction}\label{secintro}
The advances in power electronics technology bring new flexibility to power systems that facilitates system operation and control.
The flexibility can be classified into node-side flexibility (in the form of flexible power injections) and grid-side flexibility (in the form of flexible line impedances), which are enabled by shunt devices (e.g., STATCOM) and series devices (e.g., TCSC), respectively \cite{li2018grid}.
Node-side flexibility can be easily incorporated into optimal power flow (OPF) as the flexible power injections can be handled in a similar way to generation dispatch.
However, the OPF problem becomes much more complex if flexible line impedances are considered, due to that the admittance matrix is changed from a constant to a variable.
In addition, convex relaxation techniques \cite{low2014convex1, molzahn2019survey}, which have been proved powerful in addressing the non-convexity of AC-OPF problems, does not apply to the case with flexible line impedances. This adds further difficulties to the solution method and guarantee of (global) optimality.

\indent
As a consequence, the current mainstream formulations of OPF problems with flexible line impedances are based on DC power flow, which simplifies the problem but fails to capture voltage behaviors.
Moreover, the flexible line impedances still render non-convex constraints even under DC power flow. Additional binary variables, which are less friendly than continuous variables, need to be introduced to transform those non-convex constraints into linear ones \cite{ding2015optimal, sahraei2015fast}.
This paper studies the OPF with flexible line impedances based on AC power flow, and proposes a circuit transformation that makes the problem compatible with the existing convex relaxation framework.
A key discovery is made which establishes an equivalence between a flexible-impedance line and a constant-impedance line with a correlated pair of tap-adjustable transformers placed at its two terminals.
With this circuit equivalent, the original problem can be reformulated into a semi-definite program (SDP) using convex relaxation technique, which enables use of an easy-to-implement solver.
Moreover, the relaxation exactness can be enforced by adding a small penalty to the objective, which has very minor impact on the optimality, i.e., the obtained solution is near globally optimal.

\section{OPF with flexible line impedances}
Consider a power system with the set of buses $\mV=\{1,...n\}$ and set of transmission lines $\mE\subseteq\mV \times \mV$.
The set of generator buses is denoted by $\mV_G\subseteq \mV$.
The set of flexible-impedance lines is denoted by $\mE_{\textup{flex}}\subseteq \mE$, and the set of constant-impedance lines is denoted by $\mE\backslash\mEf$.
Let $V_i\in\mbC$ denote the complex voltage at bus $i$, and $P_{Gi}, P_{Li}\in\mbR$ (or $Q_{Gi}, Q_{Li}\in\mbR$) denote the active (or reactive) power generation and load at bus $i$. Set $P_{Gi}=Q_{Gi}=0$, $\forall i\notin \mV_G$.
In addition, let $y_{ij}=y_{ji}\in\mbC$ denote the admittance of line $(i,j)\in\mE$, and $y_{io}\in\mbC$ denote the shunt component at bus $i$ that consists of line charging capacitance and reactive power compensation.
Given $P_{Li}, Q_{Li}$ at each bus, the OPF problem considering flexible line impedances, say OPF-Ybus, can be formulated as
\begin{subequations}\label{OPFtcsc}
\begin{align}
    &\textup{OPF-Ybus:}
    \min_{\substack{P_{Gi}, Q_{Gi}, V_i, k_{ij}}}~\sum\nolimits_{i\in \mV_G} f_i(P_{Gi})   \label{OPFobj}\\
    &s.t.~P_{Gi}-P_{Li} + \upj (Q_{Gi}-Q_{Li}) = y_{io}^*|V_i|^2 \notag \\
    &~~~~~~~~~~~~~~~+\sum\nolimits_{j\in\mN_i} y_{ij}^* V_i (V_i^*-V_j^*),~\forall i\in \mV \label{OPFpf}\\
         &~~~~~P_{Gi}^{\min} \leq P_{Gi} \leq P_{Gi}^{\max},~\forall i\in \mV_G \label{OPFPlim}\\
         &~~~~~Q_{Gi}^{\min} \leq Q_{Gi} \leq Q_{Gi}^{\max},~\forall i\in \mV_G \label{OPFQlim}\\
         &~~~~~V_i^{\min} \leq |V_i| \leq V_i^{\max},~\forall i\in \mV   \label{OPFVlim} \\
         &~~~~~|\textup{Re}\{y_{ij}^* V_i (V_i^*-V_j^*)\}|\leq P_{ij}^{\max},~\forall (i,j)\in \mE \label{OPFPijlim} \\
         &~~~~~y_{ij} = k_{ij}\cdot\upj b_{ij}^{\textup{rated}},~\forall (i,j)\in \mEf \label{OPFyij} \\
         &~~~~~k_{ij}^{\min}\leq k_{ij} \leq k_{ij}^{\max},~\forall (i,j)\in \mEf \label{OPFkijlim}
\end{align}
\end{subequations}
where the objective \eqref{OPFobj} represents the total generation cost that is a convex quadratic function of $P_{Gi}$;
\eqref{OPFpf} refers to AC power flow equation with the superscript $*$ denoting complex conjugate and $j\in\mN_i$ meaning that bus $i$ and bus $j$ are directly connected by a line;
\eqref{OPFPlim} and \eqref{OPFQlim} refer to generation limits;
\eqref{OPFVlim} refers to voltage limits;
\eqref{OPFPijlim} refers to active line flow limits.
The flexible-impedance line is described by \eqref{OPFyij} and \eqref{OPFkijlim} where $b_{ij}^{\textup{rated}}<0$ is the rated line susceptance, $k_{ij}$ is the tuning ratio, and $0<k_{ij}^{\min}\leq 1$ and $k_{ij}^{\max}\geq 1$ are the lower and upper limits.
Let us further comment the complexity and possible extensions of OPF-Ybus as follows.

\begin{remark}[Hardness of OPF-Ybus]
For the OPF problem with fixed line parameters, the non-convexity is mainly caused by the bilinear term $V_i V_j^*$ in complex domain (or the trilinear term $|V_i||V_j|\sin(\theta_i-\theta_j)$ in real domain) in power flow equation, which has been successfully handled by convex relaxation \cite{low2014convex1, molzahn2019survey}.
When line impedances become flexible, it gives rise to the trilinear term $y_{ij}^* V_i V_j^*$ in complex domain (or the quadrilinear term $b_{ij}|V_i||V_j|\sin(\theta_i-\theta_j)$ in real domain) in power flow equation.
Although the interior-point method is a solver for general non-convex programs such like OPF-Ybus, it is generally hard to guarantee the global optimality of the obtained solution.
In addition, the common linearization technique such as the McCormick envelope does not apply to OPF-Ybus as the envelope for quadrilinear terms is not known explicitly~\cite{costa2012relaxations}.
The next section will reveal a circuit equivalence between a flexible-impedance line and a pair of ideal transformers with correlated and adjustable tap ratios.
The convex relaxation techniques, which are originally intended for conventional OPF problems, can then apply to OPF-Ybus after the circuit transformation.
Moreover, the convex reformulation of OPF-Ybus enables us to quantify the optimality in a convenient way, which will be detailed in Remark \ref{Rem3}.
\end{remark}

\begin{remark}[Modeling Extensions]
For simplicity, only line susceptance is assumed adjustable in OPF-Ybus, which can be realized by TCSC---the most common series device.
Nevertheless, the problem formulation can be easily extended to the following cases where the proposed convexification method still applies:
\begin{itemize}
  \item The line susceptance and conductance are adjustable proportionally, i.e., $y_{ij} = k_{ij}(g_{ij}^{\textup{rated}}+\upj b_{ij}^{\textup{rated}})$. It will be seen later that this functionality can be realized by power flow routers (PFRs) \cite{lin2016optimal}.
      In this case, the analysis in the next section is still valid with the term $\upj b_{ij}^{\textup{rated}}$ replaced by $g_{ij}^{\textup{rated}}+\upj b_{ij}^{\textup{rated}}$.
  \item The line susceptance and conductance can be separately adjusted, i.e., $y_{ij} = k_{ij}^{g} g_{ij}^{\textup{rated}} + \upj k_{ij}^{b} b_{ij}^{\textup{rated}}$, where $k_{ij}^{g}$ and $k_{ij}^{b}$ may take different values. This can be realized by SSSC which is able to regulate its terminal voltage as a function of line current.
      In this case, the line can be decomposed into a flexible conductive component in series with a flexible inductive component and apply the proposed method to these two components separately.
  \item Other operational constraints, e.g., apparent line flow limits, can be included without difficulty as long as they are convex with respect to power injections and line flows.
\end{itemize}
\end{remark}

\section{Convex relaxation via transformer equivalent}
\subsection{Equivalence by a correlated pair of transformers}
\begin{figure}[!h]
  \centering
  \includegraphics[width=3.3in]{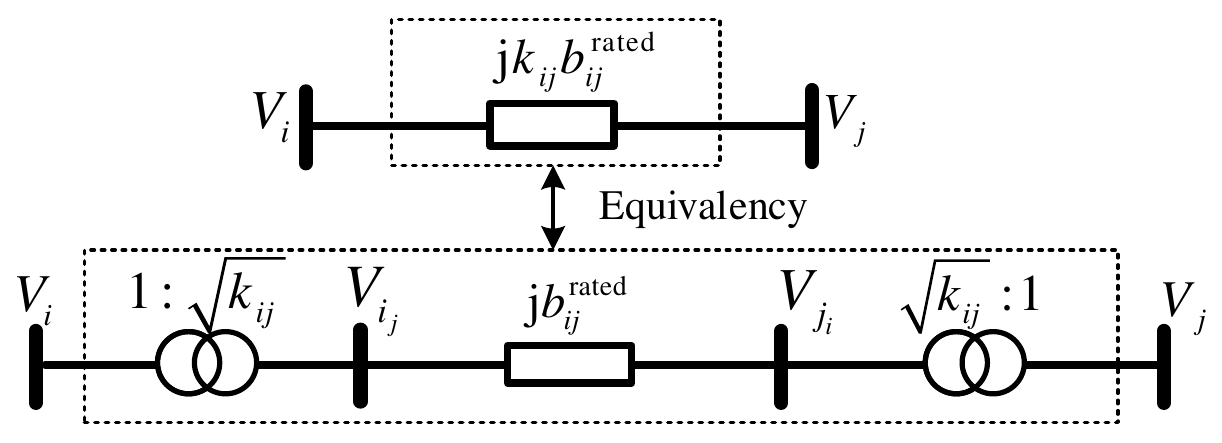}
  \caption{Equivalence between flexible line impedance and transformer pair.}
  \label{figtcscequivalent}
\end{figure}

Consider a flexible-impedance line $(i,j)$ shown in the upper circuit in Fig. \ref{figtcscequivalent}.
The complex line flows at its two sending-ends, say $S_{ij}, S_{ji}\in\mbC$, take the following form
\begin{subequations}\label{Sijtcsc}
\begin{align}
    S_{ij} &= -\upj k_{ij}b_{ij}^{\textup{rated}} V_i(V_i^*-V_j^*) \\
    S_{ji} &= -\upj k_{ij}b_{ij}^{\textup{rated}} V_j(V_j^*-V_i^*)
\end{align}
\end{subequations}
which represent the contribution of line $(i,j)$ to the power flow equation \eqref{OPFpf} of bus $i$ and bus $j$, respectively.
Clearly, the flexible-impedance line is equivalent to another device in the sense of power flow if such a device leads to the same expression of $S_{ij}, S_{ji}$ as in \eqref{Sijtcsc}.

\begin{figure}[!h]
  \centering
  \includegraphics[width=3.3in]{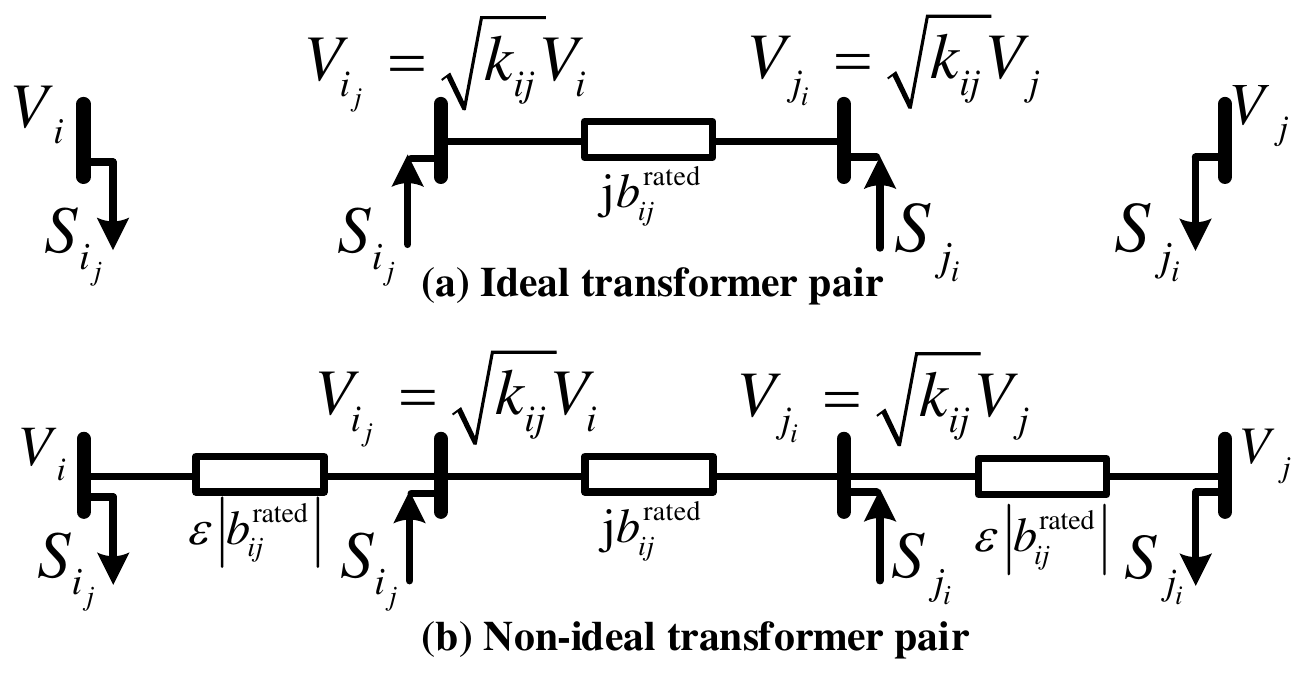}
  \caption{Power flow diagrams for transformer pair.}
  \label{figtranslossy}
\end{figure}

\indent
Then turn to the lower circuit in Fig. \ref{figtcscequivalent}, where the line has a fixed susceptance $\upj b_{ij}^{\textup{rated}}$ and a pair of ideal transformers are placed at its two terminals, say bus $i$ and bus $j$.
This pair of transformers have an identical and continuously adjustable tap ratio $\sqrt{k_{ij}}$ and the secondary-side buses adjacent to bus $i$ and bus $j$ are respectively denoted by bus $i_j$ and bus $j_i$.
The corresponding power flow diagram is shown in Fig. \ref{figtranslossy}(a), where $S_{i_j}, S_{j_i}\in\mbC$ denote the complex power transfer through the transformer pair. Observing Fig. \ref{figtcscequivalent} and Fig. \ref{figtranslossy}(a) gives
\begin{equation}\label{Sijtrans}
\begin{split}
    S_{ij} = S_{i_j} &= -\upj b_{ij}^{\textup{rated}} V_{i_j}(V_{i_j}^*-V_{j_i}^*) \\
    &= -\upj b_{ij}^{\textup{rated}} \sqrt{k_{ij}}V_i(\sqrt{k_{ij}}V_i^*-\sqrt{k_{ij}}V_j^*)\\
    &= -\upj k_{ij}b_{ij}^{\textup{rated}} V_i(V_i^*-V_j^*)
\end{split}
\end{equation}
and similarly $S_{ji} = S_{j_i} =  -\upj k_{ij}b_{ij}^{\textup{rated}} V_j(V_j^*-V_i^*)$ can be derived.
Comparing \eqref{Sijtcsc} and \eqref{Sijtrans}, this correlated pair of transformers is an equivalent model for the flexible-impedance line as they take the same effect in power flow.
Also note that installing a pair of transformers at the terminals of a line is a common realization of PFRs \cite{lin2016optimal}. So this circuit equivalence implies that flexible line impedances can be regarded as a special class of PFRs.

\indent
Further, instead of using the ideal transformer model in Fig. \ref{figtranslossy}(a) where the primary and secondary sides are fully decoupled, the non-ideal one in Fig. \ref{figtranslossy}(b) will be adopted which adds small fictitious conductance $\varepsilon|b_{ij}^{\textup{rated}}|$ and makes the primary and secondary sides still weakly coupled.
This is because the disconnection between the primary and secondary sides of ideal transformers will cause high-rank solution of the relaxed convex problem to be established in the next subsection. A detailed discussion on this issue can be found in \cite{robbins2015optimal}.
Numerically, we find that taking $\varepsilon$ in the order of $10^{-2}$ for the fictitious conductance yields satisfactory solutions with negligible error (e.g., $\varepsilon=0.04$ is adopted in the case study).

\indent
Now reformulate OPF-Ybus with all flexible-impedance lines transformed into non-ideal transformer pairs.
Since the transformer secondary-side buses $i_j, j_i$ are augmented after the circuit transformation, the power flow equation at each bus $i\in\mV$ is modified as
\begin{equation}\label{pfprimary}
\begin{split}
    &P_{Gi}-P_{Li} + \upj (Q_{Gi}-Q_{Li}) = y_{io}^*|V_i|^2 \\
    &~~~~~~~~+\sum\nolimits_{j\in\mN_i^0} y_{ij}^* V_i (V_i^*-V_j^*) \\
    &~~~~~~~~+ \sum\nolimits_{j\in\mN_i^f} S_{i_j} + \varepsilon |b_{ij}^{\textup{rated}}| V_i (V_i^*-V_{i_j}^*)
\end{split}
\end{equation}
where $j\in\mN_i^0$ means bus $i$ and bus $j$ are connected by a constant-impedance line in the original network, i.e., $(i,j)\in \mE\backslash\mEf$;
$j\in\mN_i^f$ means bus $i$ and bus $j$ are connected by a flexible-impedance line in the original network, i.e., $(i,j)\in \mEf$.
In addition, the power balance at those secondary-side buses $i_j,j_i$ takes the form below
\begin{equation}\label{pfsecondary}
\begin{split}
    S_{i_j} &= \varepsilon |b_{ij}^{\textup{rated}}| V_{i_j} (V_{i_j}^*-V_i^*) -\upj b_{ij}^{\textup{rated}} V_{i_j}(V_{i_j}^*-V_{j_i}^*)\\
    S_{j_i} &= \varepsilon |b_{ij}^{\textup{rated}}| V_{j_i} (V_{j_i}^*-V_j^*) -\upj b_{ij}^{\textup{rated}} V_{j_i}(V_{j_i}^*-V_{i_j}^*).
\end{split}
\end{equation}
Thus, OPF-Ybus can be equivalently re-expressed as
\begin{subequations}\label{OPFtcsc2}
\begin{align}
    &\textup{OPF-TP:}\min_{\substack{P_{Gi}, Q_{Gi}, S_{i_j}, V_i, V_{i_j}, k_{ij}}}~\sum\nolimits_{i\in \mV_G} f_i(P_{Gi})   \label{OPFobj2}\\
    &s.t.~\eqref{OPFVlim},~\eqref{pfprimary},~\forall i\in\mV \label{OPFpfpri}\\
         &~~~\eqref{OPFPlim},~\eqref{OPFQlim},~\forall i\in \mV_G \\
         &~~~|\textup{Re}\{y_{ij}^* V_i (V_i^*-V_j^*)\}|\leq P_{ij}^{\max},~\forall (i,j)\in \mE\backslash\mEf \label{OPFPijlim2} \\
         &~~~|\textup{Re}\{\upj b_{ij}^{\textup{rated}} V_{i_j}(V_{i_j}^*-V_{j_i}^*))\}|\leq P_{ij}^{\max},~\forall (i,j)\in \mEf  \label{OPFPijlim2f}\\
         &~~~V_{i_j}=\sqrt{k_{ij}}V_i,~V_{j_i}=\sqrt{k_{ij}}V_j,~\forall (i,j)\in \mEf \label{OPFtrans}\\
         &~~~k_{ij}^{\min}\leq k_{ij} \leq k_{ij}^{\max},~\forall (i,j)\in \mEf \label{OPFkijlim2}\\
         &~~~\eqref{pfsecondary},~\forall (i,j)\in\mEf \label{OPFpfsec}
\end{align}
\end{subequations}
where \eqref{OPFPijlim2f} is derived from \eqref{Sijtrans}, and \eqref{OPFtrans} links the primary-side and secondary-side voltages for each transformer pair.
In case the transformer tap ratios are not correlated as in \eqref{OPFtrans}, the convex relaxation has been well studied \cite{robbins2015optimal}.
The next subsection will show that the convex relaxation of OPF-TP, where the tap ratios are correlated, can be established following a similar idea.

\subsection{The relaxed SDP model}
Introduce the Hermitian matrix $\bm{W}=[W_{ij}]\in \mbC^{(n+2f)\times (n+2f)}$ where $f$ is the cardinality of $\mEf$.
Each entry of $\bm{W}$ is defined by $W_{ab}=V_a V_b^*$, where the subscripts $a,b$ are indices from the bus set $\mV$ and secondary-side buses of the transformer pairs (i.e., buses $i_j, j_i$ for $(i,j)\in\mEf$).
It is trivial that $\bm{W}$ is positive semi-definite with $\textup{rank}(\bm{W}) = 1$, from which a unique voltage profile can be recovered for the bus set $\mV$ and secondary-side buses \cite{low2014convex1}.
Then, \eqref{pfprimary} and \eqref{pfsecondary} can be rewritten in terms of $\bm{W}$ as
\begin{equation}\label{pfprimaryW}
\begin{split}
   &P_{Gi}-P_{Li} + \upj (Q_{Gi}-Q_{Li}) = y_{io}^*W_{ii} \\
   &~~~~~~~~+ \sum\nolimits_{j\in\mN_i^0} y_{ij}^* (W_{ii}-W_{ij}) \\
   &~~~~~~~~+ \sum\nolimits_{j\in\mN_i^f} S_{i_j} + \varepsilon |b_{ij}^{\textup{rated}}|(W_{ii}-W_{i i_j})
\end{split}
\end{equation}
\begin{equation}\label{pfsecondaryW}
\begin{split}
    S_{i_j} &= \varepsilon |b_{ij}^{\textup{rated}}| (W_{i_j i_j}-W_{i_j i}) -\upj b_{ij}^{\textup{rated}} (W_{i_j i_j}-W_{i_j j_i})\\
    S_{j_i} &= \varepsilon |b_{ij}^{\textup{rated}}| (W_{j_i j_i}-W_{j_i j}) -\upj b_{ij}^{\textup{rated}} (W_{j_i j_i}-W_{j_i i_j}).
\end{split}
\end{equation}
In addition, it holds the statement below regarding \eqref{OPFtrans}-\eqref{OPFkijlim2}.
\begin{prop}\label{prop1}
Constraints \eqref{OPFtrans}-\eqref{OPFkijlim2} hold if and only if the following constraints hold:
\begin{subequations}
\begin{align}
    & k_{ij}^{\min}|V_i|^2 \leq |V_{i_j}|^2 \leq k_{ij}^{\max}|V_i|^2 \label{transkijlim} \\
    & k_{ij}^{\min}|V_j|^2 \leq |V_{j_i}|^2 \leq k_{ij}^{\max}|V_j|^2 \label{transkijlim2} \\
    & V_{i_j}V_j^* = V_i V_{j_i}^*   \label{transVkij} \\
    & \textup{Im}\{V_i V_{i_j}^*\} = \textup{Im}\{V_j V_{j_i}^*\} = 0 \label{transVkij2}\\
    & \textup{Re}\{V_i V_{i_j}^*\} \geq0,~\textup{Re}\{V_j V_{j_i}^*\} \geq 0 \label{transVkij3}.
\end{align}
\end{subequations}
\end{prop}

\begin{IEEEproof}
Necessity. If \eqref{OPFtrans}-\eqref{OPFkijlim2} hold, then $V_{i_j}/V_i = V_{i_j}^*/V_i^* = V_{j_i}/V_j = V_{j_i}^*/V_j^*=\sqrt{k_{ij}}\in\mbR_{\geq0}$ that leads to \eqref{transkijlim}-\eqref{transVkij3} directly.

Sufficiency. If \eqref{transkijlim}-\eqref{transVkij3} hold, then \eqref{transVkij}-\eqref{transVkij3} give $V_{i_j}/V_i = V_{i_j}^*/V_i^* = V_{j_i}/V_j = V_{j_i}^*/V_j^*\in\mbR_{\geq0}$.
Introducing $\sqrt{k_{ij}}=V_{i_j}/V_i=V_{j_i}/V_j\in\mbR_{\geq0}$, then \eqref{transkijlim}-\eqref{transkijlim2} lead to \eqref{OPFtrans}-\eqref{OPFkijlim2}.
\end{IEEEproof}

According to the definition of matrix $\bm{W}$, \eqref{transkijlim}-\eqref{transVkij2} can be re-expressed as
\begin{subequations}\label{transW}
\begin{align}
    & k_{ij}^{\min}W_{ii}\leq W_{i_j i_j} \leq k_{ij}^{\max}W_{ii} \label{transW1} \\
    & k_{ij}^{\min}W_{jj}\leq W_{j_i j_i} \leq k_{ij}^{\max}W_{jj} \label{transW2} \\
    &W_{i_j j} = W_{i j_i}  \label{transW3} \\
    &\textup{Im}\{W_{i i_j}\} = \textup{Im}\{W_{j j_i}\} = 0 \label{transW4}\\
    &\textup{Re}\{W_{i i_j}\}\geq0,~\textup{Re}\{W_{j j_i}\}\geq0 \label{transW5}
\end{align}
\end{subequations}
where the correlation of tap ratios is fully captured by \eqref{transW3}-\eqref{transW5} using $\bm{W}$.
Thus, by Proposition \ref{prop1} and \eqref{transW}, OPF-TP is equivalent to the formulation below in terms of $\bm{W}$
\begin{subequations}\label{OPFtcscW}
\begin{align}
    &\textup{OPF-TPW:}\min_{\substack{P_{Gi}, Q_{Gi}, S_{i_j},\bm{W}}}~\sum\nolimits_{i\in \mV_G} f_i(P_{Gi})   \label{OPFobjW}\\
    &s.t.~\eqref{pfprimaryW},~\forall i\in\mV \label{OPFprimary} \\
         &~\eqref{OPFPlim},~\eqref{OPFQlim},~\forall i\in \mV_G \label{OPFPlimW}\\
         &~(V_i^{\min})^2 \leq W_{ii} \leq (V_i^{\max})^2,~\forall i\in \mV   \label{OPFVlimW} \\
         &~|\textup{Re}\{y_{ij}^* (W_{ii}-W_{ij})\}|\leq P_{ij}^{\max},~\forall (i,j)\in \mE\backslash\mEf \\
         &~|\textup{Re}\{\upj b_{ij}^{\textup{rated}} (W_{i_j i_j}-W_{i_j j_i})\}|\leq P_{ij}^{\max},~\forall (i,j)\in \mEf \\
         &~\eqref{pfsecondaryW},\eqref{transW},~\forall (i,j)\in\mEf \\
         &~\bm{W}\succeq 0 \label{OPFpsd}\\
         &~\textup{rank}(\bm{W}) = 1.  \label{OPFrank}
\end{align}
\end{subequations}
Since all constraints of OPF-TPW except \eqref{OPFrank} are convex with respect to $P_{Gi}, Q_{Gi}, S_{i_j},\bm{W}$, an SDP is obtained after relaxing \eqref{OPFrank}, where sophisticated solvers apply.
If the relaxed OPF-TPW \eqref{OPFobjW}-\eqref{OPFpsd} is exact, i.e., $\textup{rank}(\bm{W})=1$ at the optimum, then the globally optimal solution of OPF-Ybus is obtained where the tuning ratio of each line $(i,j)\in\mEf$ is given by $k_{ij}=W_{i_j i_j}/W_{ii}$.
In case of inexact relaxation, i.e., $\textup{rank}(\bm{W})\geq 2$ at the optimum, a quality and physically meaningful solution can still be obtained from the discussion below.

\begin{remark}[Relaxation Exactness]\label{Rem3}
Since OPF-Ybus has been transformed into OPF-TPW with fixed line parameters, all the existing results on enforcing relaxation exactness apply.
For instance, a rank-one solution can be enforced by replacing objective function \eqref{OPFobjW} with the penalized one below \cite{madani2015convex}
\begin{equation}\label{objpenalty}
\begin{split}
   \sum\nolimits_{i\in \mV_G} f_i(P_{Gi}) + w_q \sum\nolimits_{i\in\mV_G} Q_{Gi}
\end{split}
\end{equation}
where $w_q$ is a small positive number. The penalty term $w_q \sum\nolimits_{i\in\mV_G} Q_{Gi}$ promotes the relaxation exactness with very slight impact.
The optimality of the penalized problem can be quantified as follows. Let $P_{Gi}^\opt, P_{Gi}^\sdp, P_{Gi}^{w_q}$ be the optimal generation schemes given by the original problem OPF-Ybus, the relaxed OPF-TW (consisting of \eqref{OPFobjW}-\eqref{OPFpsd}) and penalized problem (consisting of \eqref{objpenalty}, \eqref{OPFprimary}-\eqref{OPFpsd}), respectively.
If the penalized problem is exact, it follows $\sum\nolimits_{i\in \mV_G} f_i(P_{Gi}^\sdp)\leq \sum\nolimits_{i\in \mV_G} f_i(P_{Gi}^\opt)\leq \sum\nolimits_{i\in \mV_G} f_i(P_{Gi}^{w_q})$ \cite{madani2015convex}.
Therefore, the ratio $\sum\nolimits_{i\in \mV_G} f_i(P_{Gi}^{w_q})/\sum\nolimits_{i\in \mV_G} f_i(P_{Gi}^\sdp)$, which is easy to obtain, provides an upper bound for the ratio $\sum\nolimits_{i\in \mV_G} f_i(P_{Gi}^{w_q})/\sum\nolimits_{i\in \mV_G} f_i(P_{Gi}^\opt)$ to quantify the optimality of the penalized problem.
It will be seen in the case study that the penalized problem gives a near globally optimal solution for OPF-Ybus, i.e., the ratio is very close to 1.0.
\end{remark}

\section{Case study}\label{seccase}
Take IEEE 118-bus system to test the proposed method.
The detailed system parameters can be found in the MATPOWER package~\cite{zimmerman2010matpower}.
Assume the five lines listed in Table~\ref{tab118opftcsc} have flexible impedances (the resistances of these lines are set to zero) and take $k_{ij}^{\min}=0.8$, $k_{ij}^{\max}=3.0$ for these lines.
The proportion of flexible-impedance lines is very low considering that the system has totally 186 lines, but it will be shown to significantly save generation cost.
There are 35 synchronous condensers in the system, which will have constantly zero active power generation.
To highlight the control effect, the active generation upper bounds are doubled and active line flow limits $P_{ij}^{\max}$ are all set to 200 MW.

With these settings, the relaxed OPF-TW and penalized problem (with $w_q=0.2$) are solved via CVX, and the obtained solutions are denoted with superscripts ``sdp'' and ``$w_q$'', respectively.
It turns out that $\textup{rank}(\bm{W}^\sdp) = 2$ and $\textup{rank}(\bm{W}^{w_q})=1$. So a physically meaningful solution can be obtained from the penalized problem, which is shown in Table~\ref{tab118opftcsc} (the detailed generation scheme is omitted due to page limit).
The resulting generation costs are $\sum\nolimits_{i\in \mV_G} f_i(P_{Gi}^\sdp)=132269$ and $\sum\nolimits_{i\in \mV_G} f_i(P_{Gi}^{w_q})=134555$, and the ratio is $\sum\nolimits_{i\in \mV_G} f_i(P_{Gi}^{w_q})/\sum\nolimits_{i\in \mV_G} f_i(P_{Gi}^\opt)=1.017$, which implies that the solution of the penalized problem is very close to the globally optimal solution of the original problem OPF-Ybus.

As for the conventional OPF without flexible line impedances, which corresponds to the problem with all $k_{ij}$ fixed to one, the optimal generation cost is significantly increased (see the second column of Table~\ref{tab118gencost}).
Note that some lines are congested at the conventional OPF, and the flexible line impedances effectively enlarge the feasibility region and mitigate congestion. Consequently, those cost-effective generators are better utilized under flexible line impedances, which reduces the total cost.
When an even more strict line flow limit is adopted, the benefit of flexible line impedances is more significant (see the third column of Table~\ref{tab118gencost}).

\begin{table}[!h]
\renewcommand{\arraystretch}{1.3}
  \caption{Optimal tuning for flexible line impedances}
  \label{tab118opftcsc}
  \centering
    \begin{tabular}{ccc}
    \hline\hline
    Flexible line $(i,j)$   &  $-b_{ij}^{\textup{rated}}$   & $k_{ij}$ \\
    \hline
    (23, 25)    & 12.500  & 1.256  \\
    (25, 27)    & 6.135  & 2.000  \\
    (42, 49)    & 3.096  & 1.901  \\
    (47, 69)    & 3.600  & 1.575  \\
    (100, 106)   & 4.367  & 1.270  \\
    \hline\hline
    \end{tabular}%
\end{table}%

\begin{table}[!h]
\renewcommand{\arraystretch}{1.3}
  \caption{Optimal generation cost with/without flexible line impedances}
  \label{tab118gencost}
  \centering
    \begin{tabular}{ccc}
    \hline\hline
    Model   & \tabincell{c}{Gen. cost (\$/h)\\($P_{ij}^{\max}=200$MW)}   & \tabincell{c}{Gen. cost (\$/h)\\($P_{ij}^{\max}=190$MW)} \\
    \hline
    $k_{ij}$ is flexible     &   134555  &  135891 \\
    Conventional OPF    &   138707  &  143811 \\
    Saved cost  &     4152   &   7920 \\
    \hline\hline
    \end{tabular}%
\end{table}%

\section{Conclusion}\label{secconclu}
It has been revealed that a flexible-impedance line is equivalent to a constant-impedance line linking a correlated pair of tap-adjustable transformers. With this circuit equivalent, a convex relaxation formulation for AC-OPF problems with flexible line impedances has been established.

\ifCLASSOPTIONcaptionsoff
  \newpage
\fi

{\footnotesize
\bibliographystyle{IEEEtran}
\bibliography{IEEEabrv,opftcsc}

}




\end{document}